\documentclass[12pt]{article}
 
\newcommand{\ignore}[1]{}

\usepackage{amsmath, amsthm, amsfonts, amssymb} 

\usepackage[dvips,final]{graphicx}

\usepackage{graphicx}

\usepackage{float}
\newfloat{figure}{H}{lof}
\floatname{figure}{\figurename}

\DeclareMathAlphabet{\eufrak}{U}{}{}{}  
\SetMathAlphabet\eufrak{normal}{U}{euf}{m}{n}
\SetMathAlphabet\eufrak{bold}{U}{euf}{b}{n}

\usepackage{amsmath, amsthm, amsfonts, amssymb}

\numberwithin{equation}{section}

\def\real{{\mathord{{\rm I\kern-2.8pt R}}}}        
\def\inte{{\mathord{{\rm I\kern-2.8pt N}}}}
\def\PP{{\mathord{{\rm I\kern-2.8pt P}}}}

\def\real{{\mathord{\mathbb R}}}

\def\inte{{\mathord{\mathbb N}}}

\def\R{\right}
\def\L{\left}

\newcommand{\disp}{\displaystyle}

\def\P{\mathbb{P}}
\def\E{\mathop{\hbox{\rm I\kern-0.20em E}}\nolimits}

\def\Q{\mathbb{Q}}

\newtheorem{prop}{Proposition}[section]

\newtheorem{lemma}[prop]{Lemma}
\newtheorem{definition}[prop]{Definition}

\newtheorem{theorem}[prop]{Theorem}
\newtheorem{remark}[prop]{Remark}

\title{ 
\Huge Estimation of quadratic variation for two-parameter diffusions
} 
 
\author{
\Large 
Anthony R\'eveillac\footnote{anthony.reveillac@univ-lr.fr}  
\\ 
Laboratoire de Math\'ematiques
\\ 
Universit\'e de La Rochelle 
\\ 
Avenue Michel Cr\'epeau 
\\ 
17042 La Rochelle Cedex 
\\ 
France 
}

\allowdisplaybreaks

\begin{document}
\hyphenation{func-tio-nals} 
\hyphenation{para-me-ter}
\maketitle
\begin{abstract}
\noindent 
In this paper we give a central limit theorem for the weighted quadratic variations process of a two-parameter Brownian motion. As an application, we show that the discretized quadratic variations $\sum_{i=1}^{[n s]} \sum_{j=1}^{[n t]} \left\vert \Delta_{i,j} Y \right\vert^2$ of a two-parameter diffusion $Y=(Y_{(s,t)})_{(s,t)\in[0,1]^2}$ observed on a regular grid $G_n$ is an asymptotically normal estimator of the quadratic variation of $Y$ as $n$ goes to infinity.                                                                                               
\end{abstract} 
 
\normalsize

\vspace{0.5cm}



\normalsize

\baselineskip0.7cm

\section{Introduction}
\label{section:introduction}

Many statistical properties of stochastic processes can be deduced from their weighted $p$-power variations. For a one parameter process $(Z_t)_{t\in[0,1]}$ observed at regular times $\{i/n, \; 0\leq i \leq n\}$, this quantity is defined as,
$$ \sum_{i=1}^n f\L(Z_{\frac{i-1}{n}}\R) \; \L( \Delta_i Z \R)^p, \quad\quad \Delta_i Z:=Z_{\frac{i}{n}}-Z_{\frac{i-1}{n}}.$$ 
For example, the study of the power variations process has been used by Barndorff-Nielsen and Shephard in \cite{Barndorff-NielsenShephard2003,Barndorff-NielsenShephard2004} to solve some financial econometric problems (application to econometrics are also given in \cite{Barndorff-NielsenGraversenJacodShephard}). These theoretical results were also used in several fields of application such as the estimation of the integrated volatility (see for example \cite{AitSahaliaJacodEstimators} and references therein), testing for jumps of a process observed at discrete times like for example in \cite{AitSahaliaJacodTestJumps}.\\\\
In this paper we give a central limit theorem for the weighted quadratic variations process of a two-parameter Brownian motion. More precisely consider a two-parameter Brownian motion $W=(W_{(s,t)})_{(s,t)\in[0,1]^2}$ and a deterministic and regular enough function $f:\real \to \real$, we show that
\begin{equation}
\label{eq:CLT}
n \sum_{i=1}^{[n \cdot]} \sum_{j=1}^{[n \bullet]} f\L(W_{\L(\frac{i-1}{n},\frac{j-1}{n}\R)}\R) \L( \L\vert \Delta_{i,j} W \R\vert^2 - \frac{1}{n^2} \R) \underset{n\to \infty}{\overset{law(\mathcal{S})}{\longrightarrow}}  \sqrt{2} \int_{[0,\cdot]\times[0,\bullet]}  f\L(W_{(u,v)}\R) \; dB_{(u,v)},
\end{equation}
where $B$ is a two-parameter Brownian motion independent of $W$ and $\Delta_{i,j} W$ denotes the increment of the process $W$ on the subset $\Delta_{i,j}:=\left[\frac{i-1}{n},\frac{i}{n}\right]\times \left[\frac{j-1}{n},\frac{j}{n}\right]$ of $[0,1]^2$ defined by
\begin{equation}
\label{eq:increments}
\Delta_{i,j} W:=W_{\left(\frac{i-1}{n},\frac{j-1}{n}\right)}+W_{\left(\frac{i}{n},\frac{j}{n}\right)}-W_{\left(\frac{i-1}{n},\frac{j}{n}\right)}-W_{\left(\frac{i}{n},\frac{j-1}{n}\right)}.
\end{equation}
The notation $law(\mathcal{S})$ used above in (\ref{eq:CLT}) means that the convergence is in the sense of stable convergence in law in the two-parameter Skorohod space. Furthermore we stress that the limiting process is defined on an extension of the considered probability basis. Note also that usual techniques of proof used in the one-parameter setting are no longer suitable to the two-parameter case. For example the It\^o formula for two-parameter diffusion processes cannot be applied as in the one-parameter setting due to the presence of an additional term. Consequently we chose to replace the usual stochastic calculus by the Malliavin calculus which is valid in general Gaussian context.\\\\
\noindent
As an application we deduce a central limit theorem for the quadratic variations process of a two-parameter diffusion $Y=(Y_{(s,t)})_{(s,t) \in [0,1]^2}$ observed on a regular grid $G_n$ which allows us to construct an asymptotically normal consistent estimator of the quadratic variation of $Y$. More precisely, consider a two-parameter stochastic process $(Y_{(s,t)})_{(s,t) \in [0,1]^2}$ defined by,  
\begin{equation}
\label{eq:diffusion} 
Y_{(s,t)}:= \int_{[0,s] \times [0,t]} \sigma\left(W_{(u,v)} \right) \,dW_{(u,v)} + \int_{[0,s]\times[0,t]} M_{(u,v)} \; du dv, \; (s,t) \in [0,1]^2,
\end{equation}
where $(W_{(s,t)})_{(s,t) \in [0,1]^2}$ is a two-parameter Brownian motion, $\sigma:\real \to \real$ is a bounded sufficiently smooth deterministic function and $(M_{(s,t)})_{(s,t) \in[0,1]^2}$ is a continuous adapted process. Assume that $(Y_{(s,t)})_{(s,t) \in [0,1]^2}$ is observed on the regular grid 
$$G_n:=\{(i/n,j/n), \; 0 \leq i,j \leq n\},$$ 
let 
\begin{equation}
\label{eq:sumofquadratic}
V_{(s,t)}^n:=\sum_{i=1}^{[n s]} \sum_{j=1}^{[n t]} \left\vert\Delta_{i,j} Y\right\vert^2, \quad (s,t)\in[0,1]^2, \; n \geq 1,
\end{equation}
and 
\begin{equation}
\label{eq:integratedquadraticvariation}
C_{(s,t)}:=\int_{[0,s]\times[0,t]} \sigma^2\left(W_{(u,v)} \right)\; du dv, \quad (s,t) \in [0,1]^2.
\end{equation}
Using (\ref{eq:CLT}) we show in Lemma \ref{lemma:convergenceofquadratics} that,
\begin{eqnarray}
\label{eq:CLT2}
&&\disp{n \left( \sum_{i=1}^{[n \cdot]} \sum_{j=1}^{[n \bullet]} \left\vert \Delta_{i,j} Y \right\vert^2 - \int_{[0,\cdot] \times [0,\bullet]} \sigma^2\L(W_{(u,v)}\R) \; du dv \right)}\\\nonumber 
&\underset{n\to \infty}{\overset{law(\mathcal{S})}{\longrightarrow}}&\disp{\sqrt{2} \int_{[0,\cdot]\times[0,\bullet]}  \sigma^2\L(W_{(u,v)}\R) \; dB_{(u,v)},}
\end{eqnarray}
which is then used to prove that the consistent estimator $V^n$ of $C$ is asymptotically normal (\textit{cf.} Proposition \ref{proposition:AsymptoticNormality}).\\\\
Similar results have been recently established in the one-parameter setting \cite{AitSahaliaJacodEstimators,GradinaruNourdin,Jacod1,Jacod2}, let us mention some of them. Consider a one-parameter semimartingale $(Z_t)_{t\in[0,1]}$ observed  at regular times $\{i/n, \; 0 \leq i \leq n\}$ with
\begin{equation}
\label{equation:dim1}
Z_t=z_0 + \int_0^t \sigma\left(B_s\right) \; dB_s +  \int_0^t b\left(B_s\right) \; ds, \quad t \in [0,1]
\end{equation}
where $(B_t)_{t \in [0,1]}$ is a standard Brownian motion and $\sigma:\real \to \real$ and $b:\real \to \real$ are sufficiently regular deterministic functions.\\
Gradinaru and Nourdin have shown in \cite{GradinaruNourdin} that, 
\begin{equation}
\label{equation:GradinaruNourdin}
\sqrt{n} \left( \sum_{i=1}^{[n \cdot]} \left\vert \Delta_i Z \right\vert^2 - \int_0^\cdot \sigma^2(B_s) \; ds \right) \underset{n\to \infty}{\overset{law}{\longrightarrow}}  \sqrt{2} \int_0^\cdot \sigma^2(s,\beta_s^{(1)}) \; d\beta_s^{(2)},
\end{equation}
where $\beta^{(1)}$ and $\beta^{(2)}$ are two independent Brownian motions. Note that the convergence obtained in (\ref{equation:GradinaruNourdin}) hold in the Skorohod space. Property (\ref{equation:GradinaruNourdin}) has been used by Gradinaru and Nourdin in \cite{GradinaruNourdin} to construct of a goodness-of-fit test for the integrated volatility (their results are even more general since they can be applied to diffusion processes of the form (\ref{equation:dim1}) where the diffusion and the drift terms depend on the observed process $Z$ and not only on $B$).\\\\
\noindent
In \cite{Jacod1,Jacod2}, Jacod proved functional limit theorems similar to (\ref{equation:GradinaruNourdin}) in a larger setting than in (\ref{equation:dim1}) since he considered quite general functions of the increments and the process $(Z_t)_{t\in[0,1]}$ was supposed to belong to the class of It\^o semimartingales which contains some non-continuous processes and L\'evy processes. We refer to \cite{AitSahaliaJacodEstimators} for similar results established by A\"it Sahalia and Jacod. We also mention that Nourdin in \cite{Nourdin} and Nourdin, Nualart and Tudor in \cite{NourdinNualartTudor} have studied weighted power variations of a one-parameter fractional Brownian motion. Furthermore Nourdin and Peccati in \cite{NourdinPeccati} have investigated the asymptotic behavior of weighted $p$-power variations for the iterated Brownian motion.\\\\ 
We proceed as follows. First we recall in Section \ref{section:stochastic analysis for two-parameter processes} some elements of stochastic analysis of two-parameter processes. Actually we present some definitions concerning stochastic calculus of two-parameter processes taken from \cite{IvanoffMerzbach} and the definition of the two-parameter Skorohod space initially introduced in \cite{Neuhaus} and in \cite{Straf}. Secondly, in Section \ref{section:Non-central limit theorem} we establish the central limit theorem (Theorem \ref{theorem:StableInLawConvergence}) for the weighted quadratic variations process of the two-parameter Brownian motion briefly presented in (\ref{eq:CLT}). As an application we prove in Section \ref{section:Estimation} that the consistent estimator $V^n$ (\ref{eq:sumofquadratic}) of the quadratic variation $C$ (\ref{eq:integratedquadraticvariation}) is asymptotically normal (Proposition \ref{proposition:AsymptoticNormality}). Finally we present in an appendix (Section \ref{section:appendix}) some background on set-indexed processes, extension of probability bases and on the Malliavin calculus for the two-parameter Brownian motion which are used in Sections \ref{section:Non-central limit theorem} and \ref{section:Estimation}.
 
\section{Stochastic analysis of two-parameter processes}
\label{section:stochastic analysis for two-parameter processes}

In this section we recall some definitions of two-parameter stochastic analysis which will be used in Sections \ref{section:Non-central limit theorem} and \ref{section:Estimation} and we present the two-parameter Skorohod space introduced in \cite{Neuhaus} and \cite{Straf}. 
\subsubsection*{Some elements of two-parameter stochastic calculus}

Let $(\Omega,\mathcal{F},(\mathcal{F}_z)_{z \in [0,1]^2},\P)$ be a filtered probability space. \\
We denote the partial order relation $\preceq$ on $[0,1]^2$ defined by,
$$ z' \preceq z \Leftrightarrow (s'\leq s \textrm{ and } t' \leq t), \quad z'=(s',t'), \; z=(s,t).$$ 
We also define the \textit{strong past information} filtration on $(\Omega,\mathcal{F},\P)$.
\begin{definition}
\label{Definition:Strong history filtration}
Let $z=(s,t)$ in $[0,1]^2$.
$$ \mathcal{F}_{(s,t)}^\ast:=\bigvee_{s' \leq s \textrm{ or } t' \leq t} \mathcal{F}_{(s',t')}.$$
\end{definition} 
\noindent
Until the end of this paper we assume that the following \textit{commutation condition} hold. This property is a \textit{conditional independence property} (CI in short) and corresponds to the condition (F4) of \cite{CairoliWalsh}.\\\\
\textbf{Assumption (CI):}\\
\label{Definition:(CI) hypothesis}
\noindent
\textit{The filtration $(\mathcal{F}_z)_{z \in [0,1]^2}$ is supposed to satisfy the \emph{(CI)} condition \textit{i.e.} for all $z=(s,t)$ and $z'=(s',t')$ in $[0,1]^2$ }
$$\E\left[\E\left[\cdot \vert \mathcal{F}_z\right] \vert \mathcal{F}_{[0,z] \cap [0,z']}\right]=\E\left[\cdot \vert \mathcal{F}_{(s \wedge s', t \wedge t')}\right].$$
\begin{definition}
\label{Definition:Martingale}
An $(\mathcal{F}_z)_{z \in [0,1]^2}$-adapted process $(Y_z)_{z \in [0,1]^2}$ is said to be 
\begin{itemize}
\item[i)] a martingale if for every $z$ and $z'$ in $[0,1]^2$ such that $z \preceq z'$
$$ \E[Y_{z'} \vert \mathcal{F}_{z} ]= Y_z,$$ 
\item[ii)] a strong martingale if for all $z$ and $z'$ in $[0,1]^2$ such that $z \preceq z'$
$$ \E[Y_{[z,z']} \vert \mathcal{F}_{z}^\ast]=0,$$
where $Y_{[z,z']}$ denotes the increments of $Y$ on the interval $[z,z']$.
\end{itemize}
\end{definition}
\noindent
As an example, we mention the two-parameter Brownian motion $(W_z)_{z \in [0,1]^2}$ is a strong martingale with respect to its natural filtration and a centered Gaussian process with covariance function,
$$ \E[W_{(s,t)} W_{(s',t')}]=(s \wedge s') (t \wedge t'), \quad (s,t), (s',t') \in [0,1]^2.$$ 
 
\subsubsection*{Skorohod space $\mathcal{D}([0,1]^2)$}

In the one-parameter setting, Skorohod introduced in \cite{Skorohod} four topologies known as $J_1$, $J_2$, $M_1$ and $M_2$. The topology $M_2$ is the weakest of theses topologies in the sense that convergence of a sequence $(x_n)_n$ of functions on $[0,1]$ to $x$ for $J_1$, $J_2$ or $M_1$ consists in the convergence of $(x_n)_n$ to $x$ in $M_2$ plus some additional conditions. The $M_2$ topology has been extended to the general setting of set-indexed functions by Bass and Pyke in \cite{BassPyke} whereas the $J_1$ topology has been extended to multiparameter functions by Neuhaus and Straf respectively in \cite{Neuhaus} and \cite{Straf}. The two-parameter Skorohod space (relative to $J_1$) introduced in \cite{Neuhaus} and \cite{Straf} is denoted by $\mathcal{D}([0,1]^2)$ and give an equivalent  to two-parameter functions of the notion of c\`adl\`ag functions on $[0,1]$. The set $\mathcal{D}([0,1]^2)$ can be equipped with a metric $d$ which makes it a Polish space and we denote by $\mathcal{L}_2$ the Borel $\sigma$-algebra on $(\mathcal{D}([0,1]^2),d)$. Note that as in the one-parameter setting the $J_1$ topology is stronger than the $M_2$ topology. Furthermore compact sets (relative to $J_1$) on $(\mathcal{D}([0,1]^2),d,\mathcal{L}_2)$ can be described thanks to a modulus of continuity $w$ which enables us to use techniques described in \cite{Billingsley} for one-parameter functions. We conclude this section by giving the definition of $w$. Let $f:[0,1]^2 \to \real$ be an element of $\mathcal{D}([0,1]^2)$ and $\delta>0$ we define $w(f,\delta)$ as,
\begin{equation}
\label{eq:modulusofcontinuity}
w(f,\delta):=\sup_{\|(s,t)-(s',t')\|<\delta} \vert f(s,t)-f(s',t') \vert,
\end{equation}
where $\|(s,t)-(s',t')\|:=\max\{\vert s-s' \vert;\vert t-t' \vert\}$ for $(s,t),(s',t') \in [0,1]^2$.

\section{Central limit theorem}
\label{section:Non-central limit theorem}

In this section we state and prove the functional limit theorem (Theorem \ref{theorem:StableInLawConvergence}) which will allow us to show in Section \ref{section:Estimation} that the consistent estimator $V^n$ (\ref{eq:sumofquadratic}) of the quadratic variation $C$ (\ref{eq:integratedquadraticvariation}) is asymptotically normal (Proposition \ref{proposition:AsymptoticNormality}).\\\\ 
Let $f:\real\to\real$ be a bounded and measurable deterministic function. Let a two-parameter Brownian motion $W=(W_{(s,t)})_{(s,t)\in[0,1]^2}$ defined on a probability basis $\mathcal{B}:=\L(\Omega,\mathcal{F},(\mathcal{F}_{(s,t)})_{(s,t)\in[0,1]^2},\P\R)$. Let also 
$$
\xi_{i,j}:= n \; f\L(W_{\L(\frac{i-1}{n},\frac{j-1}{n}\R)}\R) \; \L( \L\vert \Delta_{i,j}W\R\vert^2 -\frac{1}{n^2}  \R), \quad 1\leq i,j \leq n, \quad n\geq 1.
$$ 
The re-normalized weighted quadratic variations process $X^n=(X_{(s,t)}^n)_{(s,t)\in[0,1]^2}$ is defined as,
\begin{equation}
\label{eq:weightedquadrticvariations}
X_{(s,t)}^n:= \sum_{i=1}^{[n s]} \sum_{j=1}^{[n t]} \xi_{i,j}, \quad (s,t)\in [0,1]^2.
\end{equation} 
Stable convergence in law has been introduced by R\'enyi in \cite{Renyi58} and in \cite{Renyi63}. It requires some particular care, since here the limiting process $X$ is not defined on the probability basis $\mathcal{B}=(\Omega,\mathcal{F},(\mathcal{F}_z^\ast)_{z \in [0,1]^2},\P)$ on which the $X^n, \; n \geq 1$ are defined but an extension $\tilde{\mathcal{B}}:=(\tilde{\Omega},\tilde{\mathcal{F}},(\tilde{\mathcal{F}}_z)_{z \in [0,1]^2},\tilde{\P})$ of $\mathcal{B}$.  
\begin{theorem}
\label{theorem:StableInLawConvergence}
Assume that the deterministic function $f:\real\to\real$ considered above is bounded. Then $(X^n)_{n \geq 1}$ defined by (\ref{eq:weightedquadrticvariations}) converges $\mathcal{F}$-stably in law in the Skorohod space $(\mathcal{D}([0,1]^2),d,\mathcal{L}_2)$ to a non-Gaussian continuous process $X$ presented below in the proof by (\ref{eq:definitionofthelimitX}) defined on an extension of the probability basis $\mathcal{B}$.
\end{theorem}

\begin{proof}
Let us first describe the extension of $\mathcal{B}$ on which the limiting process $X$ is defined.\\
We denote by $\mathcal{B}':=(\Omega',\mathcal{F}',(\mathcal{F}_z')_{z \in [0,1]^2},\P')$ the two-parameter Wiener space, that is \\$\disp{\Omega':=\mathcal{C}^0([0,1]^2)}$ is the space of real-valued continuous functions on $[0,1]^2$ vanishing on the set $\{(s,t) \in [0,1]^2, \; s=0 \textrm{ or } t=0\}$. Then $\P'$ is the unique measure on $(\Omega',\mathcal{F}')$ under which the canonical process $(B_z)_{z \in [0,1]^2}$ on $\Omega'$ defined by, 
$$B_z(\omega'):=\omega'(z), \quad \omega' \in \Omega', \; z \in [0,1]^2,$$  
is a standard two-parameter Brownian motion.\\ 
Let the extension $\tilde{\mathcal{B}}:=(\tilde{\Omega},\tilde{\mathcal{F}},(\tilde{\mathcal{F}}_z^\ast)_{z \in [0,1]^2},\tilde{\P}) $ defined as,
$$
\left\lbrace
\begin{array}{l}
 \tilde{\Omega}:=\Omega \times \Omega^{'},\\
 \tilde{\mathcal{F}}:=\mathcal{F} \otimes \mathcal{F}^{'},\\
  (\tilde{\mathcal{F}}_z)_{z \in [0,1]^2}:=(\cap_{\rho > z} \mathcal{F}_\rho^\ast \otimes \mathcal{F}_\rho^{'})_{z \in [0,1]^2},\\
   \tilde{\P}(d\omega,dy):=\P(d\omega) \P'(dy).
\end{array}
\right.
$$
We will denote by $\E$ (respectively $\tilde{\E}$) the expectation under $\P$ (respectively under $\tilde{\P}$).\\
On $\tilde{\mathcal{B}}$ we define $(X_z)_{z \in [0,1]^2}$ as,
\begin{equation}
\label{eq:definitionofthelimitX}
X_z(\omega,\omega'):= \sqrt{2}\, \left( \int_{[0,z]} f(W_\rho(\omega)) \, dB_\rho \right)(\omega'), \quad z \in [0,1]^2.
\end{equation}
The process $X$ is a $\mathcal{F}$-progressive conditional Gaussian martingale with independent increments on $\tilde{\mathcal{B}}$, which means that $X$ is an $(\tilde{\mathcal{F}}_z)_{z \in [0,1]^2}$- adapted process such that for $\P$ almost $\omega$ in $\Omega$, $X(\omega,\cdot)$ is a Gaussian process on $\mathcal{B}'$ with covariance function
\begin{eqnarray*}
&&\disp{\E_{\P'}\left[X_{(s_1,t_1)}(\omega,\cdot) X_{(s_2,t_2)}(\omega,\cdot) \right]}\\ 
&=& \disp{2 \, \int_{[s_1 \wedge s_2,s_1 \vee s_2]\times [t_1 \wedge t_2,t_1 \vee t_2]} f^2\L(W_\rho\R)(\omega) \, d\rho, \quad (s_1,t_1), (s_2,t_2) \in [0,1]^2.}    
\end{eqnarray*}
Note that $\tilde{\mathcal{B}}$ is clearly a \textit{very good extension} of $\mathcal{B}$ in the sense of Definition \ref{Definition:Very good extension}.\\ \\
\noindent
Since $(\mathcal{D}([0,1]^2),d,\mathcal{L}_2)$ is a Polish space, by \cite[Proposition VIII.5.33]{JacodShiryaev}, $\mathcal{F}$-stable convergence in law holds if for every random variable $Z$ on $(\Omega,\mathcal{F},\P)$ the couple $(Z,X^n)_n$ converges in law. Adapting an argument presented in the proof of \cite[Theorem VIII.5.7 b)]{JacodShiryaev}, the convergence in law of a such couple $(Z,X^n)_n$ will be obtained as follows. First we give a tightness property for the sequence $(X^n)_n$ (relative to the Skorohod space $(\mathcal{D}([0,1]^2),d,\mathcal{L}_2)$) and then we make an \textquotedblleft identification of the limit\textquotedblright \textit{via} $\mathcal{F}$-stable finite-dimensional convergence in law to $X$. Recall that the latter property means that for every integer $m \geq 0$, for every continuous and bounded function $\psi:\real^{m+1}\to\real$ and every elements $z_0,\ldots,z_m$ in a dense subset of $[0,1]^2$, 
\begin{equation}
\label{equation:DefinitionStableFiniteDimensional}
\E\left[Z \psi \left(X^n(z_0),X^n(z_1),\ldots,X^n(z_m)\right)\right] \underset{n \to \infty}{\longrightarrow} \tilde{\E}\left[Z \psi \left(X(z_0),X(z_1),\ldots,X(z_m)\right)\right].
\end{equation}
The proof is decomposed in two steps. In \begin{small}\textbf{\textit{Step 1)}}\end{small} we show that $(X^n)_n$ is tight in $(\mathcal{D}([0,1]^2),d,\mathcal{L}_2)$ and in \begin{small}\textbf{\textit{Step 2)}}\end{small} we prove the $\mathcal{F}$-stable finite-dimensional convergence in law to $X$.\\\\
\begin{small}\textbf{\textit{Step 1)}}\end{small}\\
We show the sequence $(X^n)_n$ is tight in the Skorohod space $(\mathcal{D}([0,1]^2),d,\mathcal{L}_2)$.\\ 
A complete description of $(\mathcal{D}([0,1]^2),d,\mathcal{L}_2)$ can be found in \cite{Neuhaus}. In particular it is shown in \cite{Neuhaus} that the set of conditions (\ref{eq:tightness1}) and (\ref{eq:tightness2}) is necessary and sufficient for the sequence $(X^n)_n$ to be tight in $(\mathcal{D}([0,1]^2),d,\mathcal{L}_2)$,
\begin{equation}
\label{eq:tightness1}
(X_0^n)_n \textrm{ converges in distribution},
\end{equation}
\begin{equation}
\label{eq:tightness2}
\lim_{\delta \to 0} \limsup_{n \to \infty} \P[w(X^n,\delta) \geq \varepsilon]=0, \quad \varepsilon >0,
\end{equation}
where $w$ is defined in (\ref{eq:modulusofcontinuity}). Property (\ref{eq:tightness1}) is clear since for every $n \geq 1$ $X_0^n=X_0=0, \; \P$-a.s.. We will show (\ref{eq:tightness2}) using a method from \cite[p. 89]{Billingsley}.\\
Let $\varepsilon>0$, $\delta>0$ and $n \geq 1$. Let $m:=\left[\frac{n}{\delta}\right]$ and $v:=\left[\frac{n}{m}\right]$. We consider on $[0,1]^2$ the rectangles $R_{i,j}:=\left[\frac{m_{i-1}}{n},\frac{m_i}{n}\right] \times \left[\frac{m_{j-1}}{n},\frac{m_j}{n}\right], \; (i,j) \in \{1,\cdots,v\}^2$ where $m_i:=i m, \; 1 < i < v$ and $m_v=n$. With this notation the length of the shortest side of the rectangles $R_{i,j}$ is greater than $\delta$ and $v\leq 2/\delta$.
We can adapt the proof of \cite[Theorem 7.4]{Billingsley} to our case and we have,
\begin{equation}
\label{eq:(8.10)}
\P[w(X^n,\delta) \geq 3 \varepsilon] \leq \sum_{i=1}^{v} \sum_{j=1}^v \P\left[\sup_{z \in R_{i,j}} \left\vert X_z^n-X_{\left(\frac{m_{i-1}}{n},\frac{m_{j-1}}{n}\right)}^n \right\vert  \geq \varepsilon \right]. 
\end{equation}
Let us give some notations.
For $(k,j)\in \{1,\ldots,n\}^2$ let 
$$S_{k,l}:=\sum_{i=1}^k \sum_{j=1}^l f\L(W_{\L(\frac{i-1}{n},\frac{j-1}{n}\R)}\R) \left(\left\vert\Delta_{i,j} W\right\vert^2 -1/n^2\right),$$
that is $X_{(k,l)}^n=n S_{k,l}$. For $z$ in $R_{i,j}$ we write $\hat{S}_{k,l}^{i,j}:=S_{k,l}-S_{\frac{m_{i-1}}{n},\frac{m_{j-1}}{n}}$. Using these notations we can write (\ref{eq:(8.10)}) as,
\begin{eqnarray*}
\disp{\P\left[w(X^n,\delta) \geq 3 \varepsilon\right]}&\leq&\disp{\sum_{i=1}^{v} \sum_{j=1}^v \P\left[\sup_{m_{i-1} \leq k \leq m_i, \; m_{j-1} \leq l \leq m_j} \left\vert \hat{S}_{k,l}^{i,j}\right\vert \geq \frac{\varepsilon}{n} \right]} \\
\end{eqnarray*}
We will now use \cite[Section 10]{Billingsley} which provides maximal inequalities for partial sums of non-independent and non-stationary random variables. For $i,j$ fixed as above, we re-index the random variables appearing in $\hat{S}_{k,l}^{i,j}$ to obtain, 
$$ \hat{S}_{k,l}^{i,j}=\sum_{p=1}^{\eta(i,j,k,l)} \tau_p,$$
with $\tau_p$ equal to some $\xi_{\cdot,\cdot}$ divided by $n$ and $\eta(i,j,k,l)$ is an integer.\\ 
Let two integers $\alpha \leq \beta$. Since $f$ is supposed to be bounded by a non-random function, let $R:=\sup_{x \in \real} \vert f(x) \vert$ non-random. Let $K, \tilde{K}, \tilde{\tilde{K}}$ denote non-random constants.\\ 
\begin{eqnarray}
\label{eq:(10.9)}
\disp{\P\left[ \left\vert \sum_{p=\alpha+1}^\beta \tau_p \right\vert \geq \lambda \right]} &\leq& \disp{\frac{1}{\lambda^4} \E\left[ \left\vert \sum_{p=\alpha+1}^\beta \tau_p \right\vert^4 \right]}\\ \nonumber 
&=&\disp{\frac{1}{\lambda^4} \sum_{p=\alpha+1}^\beta \E[\vert  \tau_p \vert^4]} \nonumber \\
&\leq& \disp{\frac{K R^4}{\lambda^4 n^8} (\beta-\alpha)^{2 \rho},\quad \frac12<\rho<1}. 
\end{eqnarray}
Using \cite[Theorem 10.2]{Billingsley} and (\ref{eq:(10.9)}) we obtain
\begin{equation}
\label{eq:(10.10)}
\P\left[\max_{i \leq k \leq m, j \leq l \leq m} \hat{S}_{k,l}^{i,j} \geq \lambda\right]\leq \frac{K R^4 m^{4\rho}}{n^8 \lambda^4}.
\end{equation}  
Now injecting inequality (\ref{eq:(10.10)}) in (\ref{eq:(8.10)}) we have,
\begin{eqnarray*}
\disp{\P\left[w(X^n,\delta) \geq 3 \varepsilon\right]}&\leq&\disp{\frac{v^2 \tilde{K} R^4 m^{4\rho}}{n^4 \varepsilon^4}}\nonumber \\
&\leq&\disp{\frac{\tilde{\tilde{K}}m^{4\rho}}{\varepsilon^4 n^4 \delta^2}}, \quad \textrm{since } v \leq 2/\delta,\nonumber \\
&\leq&\disp{\frac{\tilde{\tilde{K}}m^{4\rho}}{\varepsilon^4 } n^{4(\rho-1)} \delta^{4\rho-2}}, \quad \textrm{since } m=[n\delta],
\end{eqnarray*}
which leads to (\ref{eq:tightness2}). 
\\\\
\begin{small}\textbf{\textit{Step 2)}}\end{small}\\
\noindent
Here we choose to consider processes $X^n$ and $X$ as set-indexed processes and we use all the notations and definitions of Subsection \ref{Some elements about set-indexed processes}. Consequently the $\mathcal{F}$-stable finite-dimensional convergence in law property (\ref{equation:DefinitionStableFiniteDimensional}) can be rewritten as follows: for every continuous and bounded function $\psi$, for every elements $C_0,\ldots,C_m$ in a dense subset of $\mathcal{A}$ (see Subsection \ref{Some elements about set-indexed processes} for definitions and notations) and for every random variable $Z$ on $(\Omega,\mathcal{F},\P)$,
\begin{equation}
\label{eq:FinteDimensionalConvergenceForSet-IndexedProcesses}
\E\left[Z \psi \left(X^n(C_0),X^n(C_1),\ldots,X^n(C_m)\right)\right] \underset{n\to\infty}{\longrightarrow} \tilde{\E}\left[Z \psi\left(X(C_0),X(C_1),\ldots,X(C_m)\right)\right].
\end{equation}
To obtain (\ref{eq:FinteDimensionalConvergenceForSet-IndexedProcesses}) we adapt \cite[Proposition 7.3.7]{IvanoffMerzbach} which allows us to replace $\mathcal{F}$-stable finite-dimensional convergence in law with $\mathcal{F}$-\textit{stable semi-functional convergence in law} that is, for every \textit{simple flow} $\varphi$ (see Definition \ref{Definition:flow}) the sequence of one-parameter processes $(X^n \circ \varphi)_n$ converges $\mathcal{F}$-stably in law to the one-parameter process $X \circ \varphi$. Let us make precise this argument.\\
Assume that \textit{stable semi-functional convergence in law} holds. We aim at showing (\ref{eq:FinteDimensionalConvergenceForSet-IndexedProcesses}). As in \cite[Proposition 7.3.7]{IvanoffMerzbach} since for every $n\geq 1$, $X^n$ is an additive process (see Definition \ref{definition:additiveprocess}) it is enough to prove (\ref{eq:FinteDimensionalConvergenceForSet-IndexedProcesses}) for elements $C_0,\cdots,C_m$ such that there exists a simple flow $\varphi$ such that for every $i\in\{1,\cdots,m\}$, $C_i=\varphi(i/m)-\varphi((i-1)/m)$. Since the sequence of one-parameter c\`adl\`ag processes $(X^n\circ \varphi)_n$ is supposed to converge $\mathcal{F}$-stably in law to $X\circ \varphi$, and since one can choose a continuous version of $X \circ \varphi$, then the projection $\pi_{(0,1/m,\ldots,1)}:\mathcal{D}([0,1])\to\real^{m+1}$ is continuous and by mapping Theorem, 
\begin{eqnarray}
\label{eq:ConvergenceFiniteDimensionalOneParameter}
&&\disp{\E\left[Z \psi\left((X^n \circ \varphi) (0),(X^n \circ \varphi) (1/m),\cdots,(X^n \circ \varphi)(1)\right)\right]}\nonumber\\
&\underset{n \to \infty}{\longrightarrow}&\disp{\int_{\Omega \times \Omega_\varphi^{'}} Z(\omega) \psi\left((X \circ \varphi)(\omega,x)(0),\cdots,(X \circ \varphi)(\omega,x)(1)\right) \P'(dx) \P(d\omega)}\nonumber\\
&=&\disp{\tilde{E}\left[ Z \psi \left( (X \circ \varphi)(0),(X \circ \varphi)(1/m),\cdots,(X \circ \varphi)(1) \right)\right]},
\end{eqnarray}
where $Z$ and $\psi$ are like in (\ref{eq:FinteDimensionalConvergenceForSet-IndexedProcesses}).\\
Consequently relation (\ref{eq:FinteDimensionalConvergenceForSet-IndexedProcesses}) holds since 
$$X^n(C_i)=(X^n \circ \varphi)(i/m)-(X^n \circ \varphi)((i-1)/m), \quad 1\leq i \leq m.$$ 
Using the argument presented above we will now prove $\mathcal{F}$-stable semi-functional convergence in law to establish $\mathcal{F}$-stable finite-dimensional convergence in law.\\\\
Let $\varphi$ be a simple flow (we write $\varphi$ as $\varphi=(\varphi_1,\varphi_2)$). We have to show that the sequence of one-parameter c\`adl\`ag processes $(X^n\circ \varphi)_n$ converges $\mathcal{F}$-stably in law to the one-parameter process $X\circ \varphi$. We give some precisions about the extension of probability basis we use. We set $\mathcal{B}_\varphi:=(\Omega,\mathcal{F},(\mathcal{F}_{\varphi(t)})_{t \in [0,1]},\P)$ and $\mathcal{B}_\varphi^{'}:=(\Omega',\mathcal{F}',(\mathcal{F}'_{\varphi(t)})_{t \in [0,1]},\P')$.  
From $\mathcal{B}_\varphi$ and $\mathcal{B}_\varphi^{'}$ we define the probability basis $\tilde{\mathcal{B}}_\varphi:=(\tilde{\Omega},\tilde{\mathcal{F}},(\tilde{\mathcal{F}}_{\varphi(t)})_{t \in [0,1]},\tilde{\P})$, with,
$$
(\tilde{\mathcal{F}}_{\varphi(t)})_{t \in [0,1]}:=(\cap_{s > t} \mathcal{F}_{\varphi(s)} \otimes \mathcal{F}'_{\varphi(s)})_{t \in [0,1]}.
$$
Let $n \geq 1$, by Lemma \ref{lemma:IM.5.1.2} and Lemma \ref{lemma:StrongMartingale} the one-parameter processes $((X^n \circ \varphi)_t)_{t \in [0,1]}$ are martingales on the probability basis $\mathcal{B}_\varphi^n:=(\Omega,\mathcal{F},(\mathcal{F}_{t}^{n,\varphi})_{t \in [0,1]},\P)$ where,
$$ \mathcal{F}_{t}^{n,\varphi}:= \mathcal{F}_{\left([ n \varphi_1(t)]n^{-1} ,[ n \varphi_2(t)] n^{-1}\right)}, \quad t \in [0,1].$$ 
We define also,
$$
\left\lbrace
\begin{array}{l}
{\mathcal{F}'}_t^{n,\varphi}:= \mathcal{F}'_{\left([ n \varphi_1(t)] n^{-1},[ n \varphi_2(t)] n^{-1}\right)}, \quad t \in [0,1],\\
\tilde{\mathcal{F}}_t^{n,\varphi}:=\cap_{s > t} \mathcal{F}_s^{n,\varphi} \otimes {\mathcal{F}'}_s^{n,\varphi},
\end{array}
\right.
$$
which lead to the following probability basis,
$$
\left\lbrace
\begin{array}{l}
{\mathcal{B}'}_\varphi^n:=(\Omega',\mathcal{F}',({\mathcal{F}'}_t^{n,\varphi})_{t \in [0,1]},\P'),\\
\tilde{\mathcal{B}}_\varphi^n:=(\tilde{\Omega},\tilde{\mathcal{F}},(\tilde{\mathcal{F}}_t^{n,\varphi})_{t \in [0,1]},\tilde{\P}).
\end{array}
\right.
$$
We will now apply \cite[Theorem IX.7.3]{JacodShiryaev} (or its triangular array formulation \cite[Theorem IX.7.28]{JacodShiryaev}). This result gives conditions insuring $\mathcal{F}$-stable finite-dimensional convergence in law for a sequence of one-parameters martingales to a continuous conditional martingale with independent increments by identifying the caracteristics of these martingales plus some additional conditions. This identification is realized in Lemma \ref{lemma:IdentificationCaracteristics} in which convergences (\ref{equation:caracteristicsII}) and (\ref{equation:caracteristicsI}) can be thought as identification of the characteristics whereas properties (\ref{eq:othogomartingale}) and (\ref{eq:othogocomposante}) ensure the $\mathcal{F}$-stable feature of the convergence. Consequently from \cite[Theorem IX.7.3]{JacodShiryaev} and Lemma \ref{lemma:IdentificationCaracteristics}, the sequence $(X^n \circ \varphi)_n$ of one-parameter martingales on $\mathcal{B}_f^n$ converges $\mathcal{F}$-stably in law to $X \circ \varphi$ on the extension $\tilde{\mathcal{B}}_\varphi$ of $\mathcal{B}_\varphi$ which ends the proof.
Note that $\tilde{\mathcal{B}}_{\varphi}$ is a very good extension of $\mathcal{B}_\varphi$ since $\tilde{\mathcal{B}}$ is a very good extension of $\mathcal{B}$. 
\end{proof}
\noindent Before turning to estimation results in Section \ref{section:Estimation} we state and prove Lemma \ref{lemma:StrongMartingale} and Lemma \ref{lemma:IdentificationCaracteristics} which were used in the proof of Theorem \ref{theorem:StableInLawConvergence}.
\begin{lemma}
\label{lemma:StrongMartingale}
We use notations of Theorem \ref{theorem:StableInLawConvergence} and of its proof. $X^n$ is a strong martingale.
\end{lemma}
\begin{proof}
Let $y$ and $x$ in $[0,1]^2$ such that $y=(y_1,y_2) \preceq x=(x_1,x_2)$. Let also, 
$$\mathcal{F}_{y}^{n,\ast}:=\mathcal{F}_{(n^{-1}[n y_1],n^{-1}[n y_2])}^\ast.$$
We have
\begin{eqnarray*}
&&\disp{\E[X^n([y,x])|\mathcal{F}_{y}^{n,\ast}]}\\
&=&\disp{n \sum_{i=[n y_1]}^{[n x_1]} \sum_{j=[n y_2]}^{[n x_2]} \E\left[ f\L(W_{(\frac{i-1}{n},\frac{j-1}{n})}\R) \left( \left\vert\Delta_{i,j} W\right\vert^2 - \frac{1}{n^2} \right) \bigg{|} \mathcal{F}_{y}^{n,\ast} \right]}\\
&=&\disp{n \sum_{i=[n y_1]}^{[n x_1]} \sum_{j=[n y_2]}^{[n x_2]} \E\left[ f\L(W_{(\frac{i-1}{n},\frac{j-1}{n})}\R) \E\left[ \left( \left\vert\Delta_{i,j} W\right\vert^2 - \frac{1}{n^2} \right) \bigg{|} \mathcal{F}_{((i-1)/n,(j-1)/n)}^\ast \right] \bigg{|} \mathcal{F}_{y}^{n,\ast} \right]}\\
&=&\disp{0.}
\end{eqnarray*}
\end{proof}

\begin{lemma}
\label{lemma:IdentificationCaracteristics}
We use notations of Theorem \ref{theorem:StableInLawConvergence} and of its proof and in particular we denote a flow $\varphi$ as $\varphi=(\varphi_1,\varphi_2).$ For every $n \geq 1$, $X^n \circ \varphi$ is a one-parameter martingale with modified second caracteristics $(0,\tilde{C}_{X^n \circ \varphi},\nu_{X^n \circ \varphi})$ on $\mathcal{B}_\varphi^n$ such that
\begin{equation} \label{equation:caracteristicsII} \nu_{X^n \circ \varphi}([0,t] \times \{ \vert x \vert > \varepsilon\})  \underset{n\to \infty}{\overset{P}{\longrightarrow}} 0, \quad \forall t \in [0,1], \; \varepsilon > 0,\end{equation} 
\begin{equation} \label{equation:caracteristicsI} \tilde{C}_{X^n \circ \varphi}(t)  \underset{n\to \infty}{\overset{P}{\longrightarrow}} 2 \int_{[0,\varphi_1(t)]\times[0,\varphi_2(t)]} f^2\L(W_\rho\R) \; d\rho, \quad \forall t \in [0,1].\end{equation}
Furthermore for every bounded martingale $\tilde{N}$ orthogonal to $W \circ \varphi$,
\begin{equation} \label{eq:othogomartingale} \langle X^n\circ \varphi, \tilde{N} \rangle_t  \underset{n\to \infty}{\overset{P}{\longrightarrow}} 0, \quad \forall t \in [0,1].\end{equation} 
We also have,
\begin{equation} \label{eq:othogocomposante} \langle X^n\circ \varphi, W \circ \varphi \rangle_t  \underset{n\to \infty}{\overset{P}{\longrightarrow}} 0, \quad \forall t \in [0,1].\end{equation}
\end{lemma}
\begin{proof}
Let $n \geq 1$. \\\\
\begin{small}\textbf{\textit{Proof of (\ref{equation:caracteristicsII}):}}\end{small}\\
In the following $\mu_{X^n \circ \varphi}$ denote the jump measure of $X^n \circ \varphi$. For $\omega \in \Omega$,
\begin{eqnarray*}
\disp{\mu_{X^n \circ \varphi}(\omega, dt, dx)}&=&\disp{\sum_{k=1}^{n } \big[ \delta_{\{\varphi_1^{-1}(k n)\}\cap (\varphi_2^{-1}(\{l n, \; 1\leq l \leq n  \}))^c,H^{1,k,n}(\omega)}(dt,dx)}\\
&+&\disp{\delta_{\{\varphi_2^{-1}(k n)\}\cap (\varphi_2^{-1}(\{l n, \; 1 \leq l \leq n\}))^c,H^{2,k,n}(\omega)}(dt,dx)\big]}\\ 
&+&\disp{ \sum_{k=1}^{n}  \sum_{l=1}^{n} \delta_{\varphi_1^{-1}(k/n),\varphi_2^{-1}(l/n),H^{3,k,l,n}(\omega)}(dt,dt,dx),}
\end{eqnarray*}
$$
\left\lbrace
\begin{array}{l}
H^{1,k,n}(\omega)=\sum_{j=1}^{[n \varphi_2(\varphi_1^{-1}(k/n))]} \xi_{k,j}(\omega), \\
H^{2,l,n}(\omega)=\sum_{i=1}^{[n \varphi_1(\varphi_2^{-1}(l/n))]} \xi_{i,l}(\omega),\\
H^{3,k,l,n}(\omega)=H^{1,k,n}+H^{2,l,n}+\xi_{k,l}.
\end{array}
\right.
$$
We denote by $\nu_{X^n \circ \varphi}$ the compensator of the measure $\mu_{X^n \circ \varphi}$. Let $A$ a Borel set in $\real$, we have
\begin{eqnarray*}
\disp{\nu_{X^n \circ \varphi}(\omega,[0,t] \times A)}&=&\disp{\sum_{k=1}^{[n \varphi_1(t)]} \E\left[\textbf{1}_{H^{1,k,n} \in A}|\mathcal{F}_{(k-1)/n,([n \varphi_2(t)]-1)/n} \right]}\\ 
&+& \disp{\sum_{l=1}^{[n \varphi_2(t)]} \E\left[\textbf{1}_{H^{2,l,n} \in A}|\mathcal{F}_{(([n \varphi_1(t)]-1)/n,(l-1)/n)} \right]}\\
&+&\disp{ \sum_{k=1}^{[n \varphi_1(t)]} \sum_{l=1}^{[n \varphi_2(t)]} \E\left[\textbf{1}_{H^{3,k,l,n} \in A}|\mathcal{F}_{([n \varphi_1(t)]-1)/n,([n \varphi_2(t)]-1)/n} \right]}.
\end{eqnarray*}
Let $\varepsilon>0$ and $k,l,n,t$ as above. Since $f$ is assumed to be bounded by a non-random constant let $R:=\sup_{x \in \real} \vert f(x) \vert$. Denote by $C$ and $\hat{C}$ some constants.  
\begin{eqnarray*}
\disp{\P\left[ \vert H^{1,k,n} \vert > \varepsilon \vert \mathcal{F}_{((k-1)/n,[n \varphi_2(t)]/n)} \right]}&\leq&\disp{\frac{1}{\varepsilon^4} \E\left[ \left\vert \sum_{j=1}^{[n \varphi_2(t)]} \xi_{k,j} \right\vert^4 \Bigg\vert \mathcal{F}_{\left(\frac{k-1}{n},\frac{[n \varphi_2(t)]-1}{n}\right)} \right]}\\
&\overset{(\ast)}{\leq}&\disp{\frac{1}{\varepsilon^4} \sum_{j=1}^{[n \varphi_2(t)]} \E\left[\left\vert \xi_{k,j} \right\vert^4 \bigg\vert \mathcal{F}_{\left(\frac{k-1}{n},\frac{[n \varphi_2(t)]-1}{n}\right)}\right]}\\
&\leq&\disp{\frac{R^8 n^4}{\varepsilon^4} \sum_{j=1}^{[n \varphi_2(t)]} \E\left[\left\vert \left(\left\vert\Delta_{k,j} W\right\vert^2-1/n^2\right) \right\vert^4\right]}\\
&\leq&\disp{\frac{R^8 C}{\varepsilon^4 n^3}}.
\end{eqnarray*}
Before giving details about the inequality $(\ast)$ we deduce of the preceding inequalities that,
\begin{eqnarray*}
&&\disp{\P\left[ \vert H^{3,k,l,n} \vert > \varepsilon \bigg\vert \mathcal{F}_{\left(\frac{[n \varphi_1(t)]-1}{n},\frac{[n \varphi_2(t)]-1}{n}\right)} \right] }\\
&\leq&\disp{\P\left[ \vert H^{1,k,n} \vert > \frac{\varepsilon}{3} \bigg\vert \mathcal{F}_{\left(\frac{[n \varphi_1(t)]-1}{n},\frac{[n \varphi_2(t)]-1}{n}\right)} \right] + \P\left[ \vert H^{2,l,n} \vert > \frac{\varepsilon}{3} \bigg\vert \mathcal{F}_{\left(\frac{[n \varphi_1(t)]-1}{n},\frac{[n \varphi_2(t)]-1}{n}\right)} \right]}\\ 
&+&\disp{\P\left[ \vert \xi_{k,l} \vert > \frac{\varepsilon}{3} \bigg\vert \mathcal{F}_{\left(\frac{[n \varphi_1(t)]-1}{n},\frac{[n \varphi_2(t)]-1}{n}\right)} \right]}\\
&\leq&\disp{\frac{\hat{C} R^8}{\varepsilon^4 n^3}}.
\end{eqnarray*}
Which leads to (\ref{equation:caracteristicsII}). Now we give some details about inequality $(\ast)$.\\
Let $1 \leq l_2,l_3,l_4 < l_1 \leq [n \varphi_2(t)]-1$. We have
\begin{eqnarray*}
&&\disp{\E\left[\xi_{k,l_1} \ldots \xi_{k,l_4}\big\vert \mathcal{F}_{\left(\frac{k-1}{n},\frac{[n \varphi_2(t)]-1}{n}\right)} \right]}\\
&=&\disp{\E\left[\E\left[\xi_{k,l_1} \ldots \xi_{k,l_4} \big\vert \mathcal{F}_{\left(\frac{k-1}{n},\frac{[n \varphi_2(t)]-1}{n}\right)} \vee \mathcal{F}_{\left(\frac{k}{n},\frac{l_1-1}{n}\right)}\right] \bigg\vert \mathcal{F}_{\left(\frac{k-1}{n},\frac{[n \varphi_2(t)]-1}{n}\right)} \right]}\\
&=&\disp{n \E\bigg[ \xi_{k,l_2} \ldots \xi_{k,l_4} f\L(W_{\left(\frac{k-1}{n},\frac{l_1-1}{n}\right)}\R) \E\bigg[ \left(\left\vert\Delta_{k,l_1} W\right\vert^2-1/n^2\right)}\\ 
&&\disp{\bigg\vert \mathcal{F}_{\left(\frac{k-1}{n},\frac{[n \varphi_2(t)]-1}{n}\right)} \vee \mathcal{F}_{\left(\frac{k}{n},\frac{l_1-1}{n}\right)}\bigg] \bigg\vert \mathcal{F}_{\left(\frac{k-1}{n},\frac{[n \varphi_2(t)]-1}{n}\right)} \bigg]}\\
&=&\disp{n \E\left[ \xi_{k,l_2} \ldots \xi_{k,l_4} f\L(W_{\left(\frac{k-1}{n},\frac{l_1-1}{n}\right)}\R) \E\left[ \left(\left\vert\Delta_{k,l_1} W\right\vert^2-1/n^2\right) \right] \bigg\vert \mathcal{F}_{\left(\frac{k-1}{n},\frac{[n \varphi_2(t)]-1}{n}\right)} \right]}\\
&=& 0.
\end{eqnarray*}
\begin{small}\textbf{\textit{Proof of (\ref{equation:caracteristicsI}):}}\end{small}\\
$\tilde{C}_{X^n \circ \varphi}(t)=\langle X^n \circ \varphi, X^n \circ \varphi \rangle$ is the compensator of $[X^n \circ \varphi, X^n \circ \varphi]$ with respect to $(\mathcal{F}_{t}^{n,\varphi})_{t \in [0,1]}$ (see for example \cite[Proof of Proposition II.2.17 b)]{JacodShiryaev}) and we have (\textit{cf.} \cite[(I.4.53)]{JacodShiryaev}),
$$ [X^n \circ \varphi,X^n \circ \varphi]_t = \sum_{0 \leq s \leq t} (X^n \circ \varphi)(s) - (X^n \circ \varphi)(s_-)=\sum_{i=1}^{[n \varphi_1(t)]} \sum_{j=1}^{[n \varphi_2(t)]} \xi_{i,j}^2,\quad  t \in [0,1].$$
Consequently, 
\begin{eqnarray*}
\disp{\tilde{C}_{X^n \circ \varphi}(t)}&=&\disp{\sum_{k=1}^{[n \varphi_1(t)]} \E\left[(H^{1,k,n})^2 \vert \mathcal{F}_{(k-1)/n,([n \varphi_2(t)]-1)/n} \right]}\\ 
&+& \disp{\sum_{l=1}^{[n \varphi_2(t)]} \E\left[(H^{2,l,n})^2 \vert \mathcal{F}_{(([n \varphi_1(t)]-1)/n,(l-1)/n)} \right]}\\
&+&\disp{ \sum_{k=1}^{[n \varphi_1(t)]} \sum_{l=1}^{[n \varphi_2(t)]} \E\left[(H^{3,k,l,n})^2 \vert \mathcal{F}_{([n \varphi_1(t)]-1)/n,([n \varphi_2(t)]-1)/n} \right]}.
\end{eqnarray*}
We can show this sum is equal to
\begin{equation}
\label{equation:caracteristicsIIbis}
\tilde{C}_{X^n \circ \varphi}(t) = \frac{2}{n^2} \sum_{i=1}^{[n \varphi_1(t)]} \sum_{j=1}^{[n \varphi_2(t)]} f^2\L(W_{\L(\frac{k-1}{n},\frac{j-1}{n}\R)}\R), \quad t \in [0,1],
\end{equation}
since terms of the form $\disp{\E\left[\xi_{k,j} \xi_{k,l} \bigg\vert \mathcal{F}_{\left(\frac{k-1}{n},\frac{[n \varphi_2(t)]-1}{n}\right)} \right]}$ vanish for $j < l \leq [n \varphi_2(t)]$ using the same type an argument described in the proof of (\ref{equation:caracteristicsII}). Furthermore terms of the form $\disp{\E\left[\xi_{k,j}^2 \vert \mathcal{F}_{\left(\frac{k-1}{n},\frac{[n \varphi_2(t)]-1}{n}\right)} \right]}$ are given by, 
\begin{eqnarray*}
\disp{\E\left[\xi_{k,j}^2 \vert \mathcal{F}_{\left(\frac{k-1}{n},\frac{[n \varphi_2(t)]-1}{n}\right)} \right]}&=&\disp{n^2 f^2\L(W_{\left(\frac{k-1}{n},\frac{j-1}{n}\right)}\R) \E\left[ \left\vert \left\vert\Delta_{k,j}^n W\right\vert^2-1/n^2 \right\vert^2 \right]}\\
&=&\disp{\frac{2}{n^2} f^2\L(W_{\L(\frac{k-1}{n},\frac{j-1}{n}\R)}\R)}.
\end{eqnarray*}
We deduce (\ref{equation:caracteristicsI}) from (\ref{equation:caracteristicsIIbis}).\\\\
\begin{small}\textbf{\textit{Proof of (\ref{eq:othogomartingale}):}}\end{small}\\
Let $\tilde{N}$ be a martingale orthogonal to $W\circ \varphi$. Without loss of generality we can assume there exists a strong martingale $N$ on $\mathcal{B}$ orthogonal to $W$ such that $NW$ is a strong martingale and such that $\tilde{N}=N\circ \varphi$. Let $n\geq 1$ and $t$ in $[0,1]$. We have,
\begin{eqnarray*}
\disp{\langle X^n, N \rangle_t}&=&\disp{\sum_{i=1}^{[n \varphi_1(t)]} \sum_{j=1}^{[n \varphi_2(t)]} \E\left[ \xi_{i,j} \Delta_{i,j} N \big{\vert} \mathcal{F}_{\left(\frac{i-1}{n},\frac{j-1}{n}\right)} \right]}\\
&=&\disp{\sum_{i=1}^{[n \varphi_1(t)]} \sum_{j=1}^{[n \varphi_2(t)]} f\L(W_{\left(\frac{i-1}{n},\frac{j-1}{n}\right)}\R) \E\left[ \left(\left\vert\Delta_{i,j}W\right\vert^2-1/n^2\right) \Delta_{i,j} N \big{\vert} \mathcal{F}_{\left(\frac{i-1}{n},\frac{j-1}{n}\right)} \right].}
\end{eqnarray*}
Let $1\leq i,j \leq n$. We use a technique presented in \cite[Lemma 6.8]{Jacod2}. For $z$ in $[0,1]^2$ we define $U_z:=\E\left[\left\vert\Delta_{i,j}W\right\vert^2-1/n^2 \vert \mathcal{F}_z\right]$. $(U_z)_{z \succeq ((i-1)/n,(j-1)/n)}$ is a martingale for the filtration generated by $U_z-U_{\left(\frac{i-1}{n},\frac{j-1}{n}\right)}$. Using the representation $\left\vert\Delta_{i,j}W\right\vert^2-\frac{1}{n^2}=I_2\left(\textbf{1}_{\Delta_{i,j}}^{\otimes 2}\right)$ as a multiple stochastic integral (see for example \cite[Section 1.1.2]{Nualart3} or \cite{NualartZakai}) we have by \cite[Lemma 1.2.5]{Nualart3} that for $z \succeq ((i-1)/n,(j-1)/n)$, $U_z=I_2\left(\textbf{1}_{\Delta_{i,j}}^{\otimes 2} \textbf{1}_{[0,z]} \right)$. From \cite{MerzbachNualart} there exists a adapted process $(\Phi_{(s,t)})_{(s,t)\in[0,1]^2}$ such that 
$$U_{(s,t)}=U_{\left(\frac{i-1}{n},\frac{j-1}{n}\right)} + \int_{[(i-1)/n,s]\times[(j-1)/n,t]} \Phi_\rho \, dW_\rho.$$
Define $N'_z=N_z-N_{\left(\frac{i-1}{n},\frac{j-1}{n}\right)}, \; z \succeq ((i-1)/n,(j-1)/n)$. This process is orthogonal to $(U_z)_{z \geq ((i-1)/n,(j-1)/n)}$. Consequently using a characterization of orthogonal two-parameter martingales given in \cite[Proposition 1.6]{CairoliWalsh} we have that 
$$\E\left[ \Delta_{i,j} N' \Delta_{i,j}U \big{\vert} \mathcal{F}_{\left(\frac{i-1}{n},\frac{j-1}{n}\right)} \right]=0.$$
A straightforward computation gives that,
$$ \E\left[ \left(\left\vert\Delta_{i,j}W\right\vert^2-1/n^2\right) \Delta_{i,j} N \big{\vert} \mathcal{F}_{\left(\frac{i-1}{n},\frac{j-1}{n}\right)} \right]=\E\left[ \Delta_{i,j} N' \Delta_{i,j}U \big{\vert} \mathcal{F}_{\left(\frac{i-1}{n},\frac{j-1}{n}\right)} \right].$$
\begin{small}\textbf{\textit{Proof of (\ref{eq:othogocomposante}):}}\end{small}\\
Let $n\geq 1$ and $t\in[0,1]$. We have, 
\begin{eqnarray*}
\disp{\langle X^n\circ \varphi,W\circ f \rangle_t}&=&\disp{\sum_{i=1}^{[n \varphi_1(t)]} \sum_{j=1}^{[n \varphi_2(t)]} \E\left[ \xi_{i,j} \Delta_{i,j} W \big{\vert} \mathcal{F}_{\left(\frac{i-1}{n},\frac{j-1}{n} \right)} \right]}\\
&=&\disp{ \sum_{i=1}^{[n \varphi_1(t)]} \sum_{j=1}^{[n \varphi_2(t)]} f\L(W_{\left(\frac{i-1}{n},\frac{j-1}{n} \right)}\R) \E\left[ \left(\left\vert\Delta_{i,j}W\right\vert^2 -\frac{1}{n^2} \right) \Delta_{i,j} W \big{\vert} \mathcal{F}_{\left(\frac{i-1}{n},\frac{j-1}{n} \right)} \right] }\\
&=&\disp{0}. 
\end{eqnarray*}
\end{proof}

\section{Estimation of the quadratic variation and asymptotic normality of the estimator}
\label{section:Estimation}

In this section we prove an asymptotic normality property (Proposition \ref{proposition:AsymptoticNormality}) for the consistent estimator $V^n$ (see (\ref{eq:sumofquadratic})) of the quadratic variation $C$ (defined in (\ref{eq:integratedquadraticvariation})).\\\\
\noindent
Consider the following two-parameter stochastic process,
\begin{equation}
\label{eq:definitionsemimartingale}
Y_{(s,t)}:=\int_{[0,s]\times[0,t]} \sigma(W_\rho) \; dW_\rho + \int_{[0,s]\times[0,t]} M_\rho \; d\rho, \quad (s,t)\in [0,1]^2
\end{equation}
defined on a probability basis $\left(\Omega,\mathcal{F},(\mathcal{F}_{(s,t)})_{(s,t)\in [0,1]^2},\P\right).$ Until the end of this paper we assume that $(M_{(s,t)})_{(s,t) \in [0,1]^2}$ is a continuous and $(\mathcal{F}_{(s,t)})_{(s,t)\in [0,1]^2}$-adapted process. Let the following assumptions which will used in the results presented below.\\\\
\textbf{Assumption (R1):}\\
\textit{The function $\sigma(\cdot)$ and its derivatives are assumed to be bounded by non-random constants and $\sigma(\cdot)$ is assumed to be at least in $\mathcal{C}^4$. }\\\\
\textbf{Assumption (R2):}\\
\textit{The function $\sigma(\cdot)$ and its derivatives are assumed to be bounded by non-random constants and $\sigma(\cdot)$ is assumed to be at least in $\mathcal{C}^8$.}\\\\
Let us define the quantities we will study.\\
Let for $n \geq 1$,
\begin{eqnarray}
\label{eq:Yn}
\disp{Y_{(s,t)}^n}&:=&\disp{n \left( V_{(s,t)}^n-C_{(s,t)} \right)}\nonumber\\
&=&\disp{n \left( \sum_{i=1}^{[n s]} \sum_{j=1}^{[n t]} \left\vert\Delta_{i,j} Y\right\vert^2 - \int_{[0,s]\times[0,t]} \sigma^2\left(W_{(u,v)} \right)\; du dv \right).}
\end{eqnarray}
\begin{prop}
\label{prop:consistent}
The estimator $V^n$ defined in (\ref{eq:sumofquadratic}) of the quadratic variation $C$ (\ref{eq:integratedquadraticvariation}) is consistent that is for every $(s,t)$ in $[0,1]^2$,
$$ V_{(s,t)}^n \underset{n\to \infty}{\overset{\P}{\longrightarrow}} C_{(s,t)}.$$
\end{prop} 
\begin{proof}
Fix $(s,t)$ in $[0,1]^2$. We show $(V_{(s,t)}^n)_n$ converges to $C_{(s,t)}$ in $L^2(\Omega,\mathcal{F},\P)$. It is enough to show that 
$$ \sum_{i=1}^{[n s]} \sum_{j=1}^{[n t]} \sigma^2\L(W_{\left(\frac{i-1}{n},\frac{j-1}{n}\right)}\R) \left( \left\vert \Delta_{i,j}W \right\vert^2 - \frac{1}{n^2} \right)$$
tends to zero in $L^2(\Omega,\mathcal{F},\P)$ as $n$ goes to infinity. We have,
\begin{eqnarray*}
&&\disp{\E\left[ \left\vert \sum_{i=1}^{[n s]} \sum_{j=1}^{[n t]} \sigma^2\L(W_{\left(\frac{i-1}{n},\frac{j-1}{n}\right)}\R) \left( \left\vert \Delta_{i,j}W \right\vert^2 - \frac{1}{n^2} \right) \right\vert^2 \right]}\\
&=&\disp{\sum_{i=1}^{[n s]} \sum_{j=1}^{[n t]} \E\left[  \left\vert  \sigma^2\L(W_{\left(\frac{i-1}{n},\frac{j-1}{n}\right)}\R) \left( \left\vert \Delta_{i,j}W \right\vert^2 - \frac{1}{n^2} \right) \right\vert^2 \right]}\\
&=&\disp{\sum_{i=1}^{[n s]} \sum_{j=1}^{[n t]} \E\left[  \sigma^4\L(W_{\left(\frac{i-1}{n},\frac{j-1}{n}\right)}\R) \E\left[ \left\vert\left\vert \Delta_{i,j}W \right\vert^2 - \frac{1}{n^2} \right\vert^2 \bigg\vert \mathcal{F}_{\left(\frac{i-1}{n},\frac{j-1}{n}\right)} \right] \right]}\\
&=&\disp{\sum_{i=1}^{[n s]} \sum_{j=1}^{[n t]} \frac{1}{n^4} \E\left[  \sigma^4\L(W_{\left(\frac{i-1}{n},\frac{j-1}{n}\right)}\R) \right] \E\left[ \left\vert\left\vert \frac{1}{n} \Delta_{i,j}W \right\vert^2 - 1 \right\vert^2 \right] }\\
&\leq&\disp{\frac{2 R^4}{n^2}},
\end{eqnarray*} 
where $R=\sup_{x \in \real}\vert  \sigma(x) \vert$.
\end{proof}
\noindent
We state and prove that the estimator $V^n$ of $C$ is asymptotically normal.
\begin{prop}[Asymptotic normality]
\label{proposition:AsymptoticNormality}
Let for $(s,t)$ in $[0,1]^2$ $S_{(s,t)}^n$ be
$$S_{(s,t)}^n:=n^2 \sum_{i=1}^{[n s]} \sum_{j=1}^{[n t]} \left\vert\Delta_{i,j}Y\right\vert^4, \quad n \geq 1.$$ 
Under \emph{\textbf{assumption (R2)}} for $(s,t)$ fixed in $(0,1]^2$ we have,
$$ \left(S_{(s,t)}^n \right)^{-\frac12} \, n \, (V_{(s,t)}^n-C_{(s,t)}) \underset{n\to \infty}{\overset{law}{\longrightarrow}} \sqrt{\frac{2}{3}} N, \quad N \sim \mathcal{N}(0,1).$$
\end{prop}
\begin{proof}
Using a localization argument, only the stochastic integral part of $Y$ gives a contribution to the limit so we assume $M=0$ in (\ref{eq:definitionsemimartingale}). The main argument of the proof is the convergence in law of $(S^n,Y^n)_n$ to $(S,X)$ where $S$ is defined by,
$$S_{(s,t)}:=3 \int_{[0,s]\times[0,t]} \sigma^4\left(W_\rho\right) \, d\rho.$$
Actually assume this convergence hold. Since $(x,y)\mapsto x^{-\frac12} y$ is continuous on $\real_+^\ast \times \real$ we have for every $(s,t)$ in $(0,1]^2$, 
$$ \left(S_{(s,t)}^n \right)^{-\frac12} \, Y_{(s,t)}^n  \underset{n\to \infty}{\overset{law}{\longrightarrow}} \sqrt{\frac23} \frac{\int_{[0,s] \times [0,t]} \sigma^2 \left(W_\rho \right) \, dB_\rho}{\left( \int_{[0,s] \times [0,t]} \sigma^4\left( W_\rho \right) \, d\rho \right)^{-\frac12}}.$$ 
Computing the characteristic function with respect to the probability measure $\tilde{\P}$ we can show that
$$\sqrt{\frac23} \frac{\int_{[0,s] \times [0,t]} \sigma^2 \left(W_\rho \right) \, dB_\rho}{\left( \int_{[0,s] \times [0,t]} \sigma^4\left( W_\rho \right) \, d\rho \right)^{-\frac12}}\overset{law}{=}\sqrt{\frac23} N, \quad N \sim \mathcal{N}(0,1).$$
We have now to show that $(S^n,Y^n)_n$ converges in law to $(S,X)$. The key point is the $\mathcal{F}$-stable convergence in law of $(Y^n)_n$ to $X$ obtained in Lemma \ref{lemma:convergenceofquadratics} stated and proved at the end of this section. Using a result of Aldous and Eagleson (presented in \cite{AldousEagleson}) concerning stable convergence in law if $(S^n)_n$ converges in $\P$-probability to $S$ then $(S^n,Y^n)_n$ converges in law to $(S,X)$ (and the convergence is even $\mathcal{F}$-stable convergence in law). \\
Let us finally show that $(S^n)_n$ converges in $L^2(\Omega,\mathcal{F},\P)$ to $S$.\\
First we show that
$$ \E\left[ \left\vert n^2 \sum_{i=1}^{[n s]} \sum_{j=1}^{[n t]} \sigma^4\L(W_{\left(\frac{i-1}{n},\frac{j-1}{n} \right)}\R) \left( \left\vert\Delta_{i,j} W \right\vert^4-\frac{3}{n^4} \right) \right\vert^2 \right] \underset{n\to \infty}{\longrightarrow} 0.$$
Actually for every $(s,t)$ in $[0,1]^2$,
\begin{eqnarray*}
&&\disp{\E\left[ \left\vert n^2 \sum_{i=1}^{[n s]} \sum_{j=1}^{[n t]} \sigma^4\L(W_{\left(\frac{i-1}{n},\frac{j-1}{n} \right)}\R) \left( \left\vert\Delta_{i,j} W \right\vert^4-\frac{3}{n^4} \right) \right\vert^2 \right]} \\
&=&\disp{n^4 \sum_{i=1}^{[n s]} \sum_{j=1}^{[n t]} \E\left[ \sigma^8\L(W_{\left(\frac{i-1}{n},\frac{j-1}{n}\right)}\R) \left\vert \left\vert\Delta_{i,j} W \right\vert^4-\frac{3}{n^4} \right\vert^2 \right]}\\
&\leq&\disp{\frac{C}{n^2}},
\end{eqnarray*}
where $C$ is a constant.\\
Using a Riemann approximation for integrals we have that
$$ 3 n^2 \sum_{i=1}^{[n s]} \sum_{j=1}^{[n t]} \frac{\sigma^4\L(W_{\left(\frac{i-1}{n},\frac{j-1}{n} \right)}\R)}{n^4}  \underset{n\to \infty}{\overset{\P}{\longrightarrow}} 3 \int_{[0,s]\times[0,t]} \sigma^4\left( W_\rho \right) \, d\rho, \quad \forall (s,t) \in [0,1]^2.$$
The proof is finished if we can show the estimate (\ref{eq:approximationbis})
\begin{equation}
\label{eq:approximationbis}
\E \left[ \left\vert n^2\sum_{i=1}^{[n s]} \sum_{j=1}^{[n t]} \left(\left\vert \Delta_{i,j} Y \right\vert^4 - \frac{\sigma^4\L(W_{\left(\frac{i-1}{n},\frac{j-1}{n} \right)}\R)}{n^4} \right) \right\vert^2 \right] \underset{n\to \infty}{\longrightarrow} 0.
\end{equation}
Using the same techniques (successive Malliavin integrations by parts, estimates of the form (\ref{eq:estimate1})) and (\ref{eq:estimate2}) we obtain (\ref{eq:approximationbis}).\\
\end{proof}
\noindent
Using classical techniques the asymptotic normality property of $V^n$ enable construction of confidence interval for $C_{(s,t)}$ for every $(s,t)$ in $[0,1]^2$.
\begin{remark}
It seems that conditions on the regularity on $\sigma(\cdot)$ imposed in Theorem \ref{theorem:StableInLawConvergence} and in Proposition \ref{proposition:AsymptoticNormality} are too strong. These type of regularities are imposed by the Malliavin calculus techniques we use. However different techniques used in the one-parameter case like for example the It\^o formula are much more difficult to use in the two-parameter setting. In contradistinction, the Malliavin calculus presents a large range of properties valid in general Gaussian context. 
\end{remark}
\noindent
In order to prove the asymptotic normality of $V^n$ in Proposition \ref{proposition:AsymptoticNormality} we have used the following lemma.
\begin{lemma}
\label{lemma:convergenceofquadratics}
Let $X$ be the process defined in (\ref{eq:definitionofthelimitX}) where $f$ in (\ref{eq:definitionofthelimitX}) is replaced by $\sigma^2$. Under \textbf{Assumption (R1)}, $(Y^n)_{n \geq 1}$ defined by (\ref{eq:Yn}) converges $\mathcal{F}$-stably in law in the Skorohod space $(\mathcal{D}([0,1]^2),d,\mathcal{L}_2)$ to the non-Gaussian continuous process $X$ defined on the extension $\tilde{\mathcal{B}}$ described in the proof of Theorem \ref{theorem:StableInLawConvergence}.
\end{lemma} 
\begin{proof}
Using a localization argument the finite variation part of $Y$ has no contribution in the limit. So we can assume that $M=0$.\\
From Theorem \ref{theorem:StableInLawConvergence} with $f=\sigma^2$, the process $(X^n)_n$ converges $\mathcal{F}$-stably in law to $X$ with,
$$X_{(s,t)}^n=n \sum_{i=1}^{[n s]} \sum_{j=1}^{[n t]} \sigma^2\L(W_{\L(\frac{i-1}{n},\frac{j-1}{n}\R)}\R) \L(\L\vert \Delta_{i,j}W \R\vert^2 -\frac{1}{n^2} \R), \quad (s,t)\in[0,1]^2.$$ 
To conclude the proof we show that $Y^n$ is equal to $X^n$ plus a term $r_n$ which become negligible when $n$ to infinity. More precisely using the notations
$$
\left\lbrace
\begin{array}{l}
\eta_{i,j}:=\left( \int_{\Delta_{i,j}} \sigma\L(W_u\R) \, dW_u \right)^2 - \sigma^2\L(W_{\left(\frac{i-1}{n},\frac{j-1}{n}\right)}\R) \left\vert \Delta_{i,j}W\right\vert^2,\\
\eta_{i,j}':=-\int_{\Delta_{i,j}} \sigma^2\L(W_\rho\R) - \sigma^2\L(W_{\left(\frac{i-1}{n},\frac{j-1}{n}\right)}\R) \, d\rho, \\
\; 1 \leq i,j \leq n, \; n \geq 1,
\end{array}
\right.
$$
$r_n$ can be decomposed as, 
\begin{equation}
\label{eq:decomposition of r_n}
r_n(s,t):=r_n^{(1)}(s,t)+r_n^{(2)}(s,t)
\end{equation}
where
$$
\left\lbrace
\begin{array}{l}
r_n^{(1)}(s,t):=\sum_{i=1}^{[n s]} \sum_{j=1}^{[n t]} \eta_{i,j} + \eta_{i,j}', \\
r_n^{(2)}(s,t):=-\int_{[n s]/n}^s \int_{[n t]/n}^t \sigma^2\L(W_\rho\R) \,d\rho.   
\end{array}
\right.
$$ 
Using a standard argument of the form \cite[Theorem 3.1]{Billingsley} as $(X^n)_n$ converges $\mathcal{F}$-stably in law to $X$ it is enough to prove that $n \sup_{(s,t)\in [0,1]^2} \vert r_n(s,t)\vert$ converges in probability to zero to obtain the $\mathcal{F}$-stable convergence in law of $(Y^n)_n$  to $X$. We use the decomposition (\ref{eq:decomposition of r_n}) and we show
\begin{equation}
\label{eq:sup r_n^(1)}
n \sup_{(s,t)\in [0,1]^2} \vert r_n(s,t)^{(1)} \vert \overset{\P}{\rightarrow} 0,
\end{equation}
\begin{equation}
\label{eq:sup r_n^(2)}
n \sup_{(s,t)\in [0,1]^2} \vert r_n(s,t)^{(2)} \vert \overset{\P}{\rightarrow} 0.
\end{equation} 
\begin{small}\textbf{\textit{Proof of (\ref{eq:sup r_n^(1)}):}}\end{small}\\
Using Burkholder's inequality for two-parameter martingales (see Remark 2 of \cite{Nualart1}) it is enough to show that 
\begin{equation}
\label{eq:sup r_n^(1)bis}
\E\left[ \left\vert \eta_{i,j}+\eta_{i,j}' \right\vert^2 \right] \leq \frac{C}{n^5},
\end{equation}
where $C$ is a constant.\\
The tool used here is the Malliavin calculus (see Appendix \ref{subsection:Malliavin calculus for two-parameter Brownian motion}) and especially the Malliavin integration by parts formula (\ref{eq:IntegrationByPartsFormula}). The main problem comes from the computation of $\E\left[ \left\vert \eta_{i,j}\right\vert^2 \right]$. First we express $\eta_{i,j}$ as,
\begin{eqnarray*}
\disp{\eta_{i,j}}&=&\disp{\left( \int_{\Delta_{i,j}} \sigma\L(W_\rho\R)-\sigma\L(W_{\left(\frac{i-1}{n},\frac{j-1}{n}\right)}\R) \, dW_\rho\right) \, \left( \int_{\Delta_{i,j}} \sigma\L(W_\rho\R)+\sigma\L(W_{\left(\frac{i-1}{n},\frac{j-1}{n}\right)}\R) \, dW_\rho\right)}\\
&=&\disp{\delta(u) \delta(v),}
\end{eqnarray*}
where 
$$
\left\lbrace
\begin{array}{l}
u_{(s,t)}:=\textbf{1}_{(s,t) \in \Delta_{i,j}} \left( \sigma\L(W_{(s,t)}\R) - \sigma\L(W_{\left(\frac{i-1}{n},\frac{j-1}{n}\right)}\right)\R)\\
v_{(s,t)}:=\textbf{1}_{(s,t) \in \Delta_{i,j}} \left( \sigma\L(W_{(s,t)}\R) + \sigma\L(W_{\left(\frac{i-1}{n},\frac{j-1}{n}\right)}\right)\R)
\end{array}
\right.,\quad (s,t) \in [0,1]^2.
$$
$\delta(u)$ denotes the Skorohod integral of the process $(u_{(s,t)})_{(s,t)\in [0,1]^2}$ which coincides since $u$ is adapted with the It\^o stochastic integral. Consequently, 
\begin{eqnarray}
\disp{\E\left[\left\vert \eta_{i,j} \right\vert^2\right]}&=&\disp{\E[\delta(u) \left(\delta(u) \, \delta(v)^2 \right)]}\nonumber\\
&=&\disp{ \E\left[ \langle u, \,D(\delta(u) \delta(v)^2)\rangle_{L^2([0,1]^2)} \right], \quad \textrm{by (\ref{eq:IntegrationByPartsFormula})} }\nonumber\\
&=&\disp{ \int_{\Delta_{i,j}} \E\left[u_{(s,t)}\, D_{(s,t)}(\delta(u)\, \delta(v)^2) \right] ds dt }\nonumber\\
&=&\disp{ \int_{\Delta_{i,j}} \E\left[u_{(s,t)} \, \delta(v)^2\, D_{(s,t)}(\delta(u)) \right] ds dt }\label{eq:detailsdecalculs}\\
&+& \disp{2 \int_{\Delta_{i,j}} \E\left[u_{(s,t)} \, \delta(u) \, \delta(v) \, D_{(s,t)}(\delta(v)) \right] ds dt,}\nonumber
\end{eqnarray} 
where the last equality is deduced from (\ref{eq:chainrule}). We can then apply successive Malliavin integration by parts to each term of the right hand of (\ref{eq:detailsdecalculs}). Using the \textquotedblleft Heisenberg commutativity relationship\textquotedblright (\ref{eq:Heisenbergcommutativityrelationship}) and estimates (\ref{eq:estimate1}) and (\ref{eq:estimate2}), we obtain 
$$\E\left[\left\vert \eta_{i,j} \right\vert^2\right] \leq \frac{C}{n^5},$$
where $C$ denotes a constant different from the one presented in (\ref{eq:sup r_n^(1)bis}).
\begin{equation}
\label{eq:estimate1}
\E\left[\left\vert W_{\rho}-W_{\left(\frac{i-1}{n},\frac{j-1}{n}\right)} \right\vert^k\right]\leq \frac{C_k}{n^k}, \quad \rho \in \Delta_{i,j}, \quad k=2,4,6,8,
\end{equation}
where the $C_k$ denotes some constants.\\
Let $\varphi:\real \to \real$ a deterministic function in $\mathcal{C}^2$,
\begin{eqnarray}
\label{eq:estimate2}
&&\disp{\varphi\left(W_\rho\right)-\varphi\left(W_{\left(\frac{i-1}{n},\frac{j-1}{n}\right)}\right)}\nonumber\\
&=&\disp{\varphi'\left(W_{\left(\frac{i-1}{n},\frac{j-1}{n}\right)}\right)  \left(W_\rho-W_{\left(\frac{i-1}{n},\frac{j-1}{n}\right)}\right)}\nonumber\\
&+&\disp{ \varphi''\left(W_{\left(\frac{i-1}{n},\frac{j-1}{n}\right)}\right)  \left(W_\rho-W_{\left(\frac{i-1}{n},\frac{j-1}{n}\right)}\right)^2+\tilde{r}^n},
\end{eqnarray}
where $\tilde{r}^n$ is a negligible term.\\ 
Finally note that the constant $C$ in (\ref{eq:sup r_n^(1)}) depends on the following deterministic bounds
$$
\left\lbrace
\begin{array}{l}
\sup_{x\in \real} \vert \sigma(x) \vert,\\
\sup_{\rho \in [0,1]^2} \sup_{(s,t)\in [0,1]^2} \vert D_\rho \sigma\L(W_{(s,t)}\R) \vert=\sup_{x \in \real} \vert \sigma'(x)\vert\\ 
\sup_{(\rho,\alpha) \in ([0,1]^2)^2} \sup_{(s,t)\in [0,1]^2} \vert D_\alpha D_\rho \sigma \L(W_{(s,t)}\R) \vert=\sup_{x \in \real} \vert \sigma''(x)\vert,
\end{array}
\right.
$$
where $D$ denote the Malliavin derivative (see (\ref{eq:MalliavinDerivative})).\\
\begin{small}\textbf{\textit{Proof of (\ref{eq:sup r_n^(2)}):}}\end{small}\\
This convergence result is obtained from Markov inequality and by standard arguments in numerical analysis.
\end{proof}

\section{Appendix}
\label{section:appendix}
In this section we present some definitions and results used in Sections \ref{section:Non-central limit theorem} and \ref{section:Estimation}. First we provide some background on set-indexed processes and on extensions of probability bases. Finally we briefly present the Malliavin calculus for two-parameter Brownian motion, including the Malliavin integration by parts formula which has used to obtain the estimates of the previous sections.
\subsection{Set-indexed processes}
\label{Some elements about set-indexed processes}

In the following definition it will be convenient to think of two-parameter processes on $[0,1]^2$ as set-indexed processes on $ \disp{\mathcal{A}:=\{[0,z], \; z \in [0,1]^2\}}$ where $Y_{[0,z]}:=Y_z$ and $\mathcal{F}_{[0,z]}:=\mathcal{F}_{z}$, $z$ in $[0,1]^2$. We will use indifferently one of these two points of view. 
\begin{definition}(\cite[Definition 7.3.1]{IvanoffMerzbach})
\label{Definition:flow}
Let $\mathcal{A}(u)$ be the set of all finite unions of sets from $\mathcal{A}$. A simple flow $\varphi$ is an application $\varphi:[0,1] \to \mathcal{A}(u)$ which satisfies the following properties,
\begin{itemize}
\item[i)] $\varphi$ is increasing,
\item[ii)] $\varphi$ is continuous,
\item[iii)] $\varphi(0)=\emptyset$,
\item[iv)] for some $k$ in $\inte$ there exists some increasing functions $\varphi_i:[0,1] \to \mathcal{A}$ , $i=1,\ldots, k$ such that for $\frac{i-1}{k} \leq s \leq \frac{i}{k}$, $1 \leq i \leq k$ 
$$ \varphi(s)=\varphi\left(\frac{i-1}{k} \right)\cup \varphi_i(s).$$ 
\end{itemize} 
\end{definition}
\begin{definition}(\cite[Definition 1.4.5]{IvanoffMerzbach})
\label{definition:additiveprocess}
Let $\mathcal{C}$ denote the set of elements of the form $A\setminus B$ with $A$ in $\mathcal{A}$ and $B$ in $\mathcal{A}(u)$ (the set of all finite unions of elements of $\mathcal{A}$).\\
A process $(X_z)_{z \in [0,1]^2}$ identified with $(X_A)_{A \in \mathcal{A}}$ is called additive if for every $C,C_1,C_2$ in $\mathcal{C}$ with $C=C_1\cup C_2$ and $C_1 \cap C_2 = \emptyset$ we have,
$$X_C=X_{C_1}+X_{C_2}.$$
\end{definition}
\noindent
We recall a particular case of \cite[Lemma 5.1.2]{IvanoffMerzbach}.
\begin{lemma}
\label{lemma:IM.5.1.2}
Let $X$ be a strong martingale and let $\varphi$ be a $\mathcal{A}(u)$-valued simple flow defined on $[0,1]$. Then the one-parameter process $(X_{\varphi(s)})_{s \in [0,1]}$ is a martingale with respect to the filtration $(\mathcal{F}_{\varphi(s)})_{s \in [0,1]}$. 
\end{lemma}

\subsection{Extension of a probability basis}
\label{subsection:extension of a probability basis}

We give the definition of a \textit{very good extension} of a probability basis which has been introduced in \cite{JacodMemin}. The following definition is an adaptation of \cite[Definition II.7.1]{JacodShiryaev}.
\begin{definition}
\label{Definition:Very good extension}
A probability basis $\tilde{\mathcal{B}}=(\tilde{\Omega},\tilde{\mathcal{F}},(\tilde{\mathcal{F}}_z)_{z \in [0,1]^2},\tilde{\P})$ is an extension of the probability basis $\mathcal{B}=(\Omega,\mathcal{F},(\mathcal{F}_z)_{z \in [0,1]^2},\P)$ if there exists an auxiliary probability basis $\mathcal{B}'=(\Omega',\mathcal{F}',(\mathcal{F}'_z)_{z \in [0,1]^2},\P')$ such that 
$$
\left\lbrace
\begin{array}{l}
\tilde{\Omega}=\Omega \times \Omega',\\
\tilde{\mathcal{F}}=\mathcal{F} \otimes \mathcal{F}',\\
\tilde{\mathcal{F}}_z=\cap_{\rho \preceq z} \mathcal{F}_\rho \otimes \tilde{\mathcal{F}}'_z, \; z \in [0,1]^2,\\
\tilde{\P}(d\omega,d\omega')=\P(d\omega) \Q_\omega(d\omega'),
\end{array}
\right.
$$
where $\Q_\omega(d\omega')$ is a transition probability from $(\Omega,\mathcal{F})$ into $(\Omega',\mathcal{F}')$.\\
This extension is called \textit{very good} is for every $z$ in $[0,1]^2$ and for all element $A'$ in $\mathcal{F}'_z$ $\omega \mapsto \Q_\omega(A')$ is equal $\P$-a.s. to an $\mathcal{F}_z$-measurable random variable.
\end{definition}
\noindent\textbf{Notation:}\\
\textit{If $(X_{z})_{z \in [0,1]^2}$ is a stochastic process on $\mathcal{B}$ we will denote by $(X_{z})_{z \in [0,1]^2}$ again the stochastic process defined on $\tilde{\mathcal{B}}$ by} 
$$ X_{z}(\omega,\omega'):=X_{z}(\omega), \quad z \in [0,1]^2, \; (\omega,\omega')\in \tilde{\Omega}.$$
\begin{lemma}
\label{lemma:verygoodextension}
As in \cite[Lemma II.7.3]{JacodShiryaev}, under \textbf{Assumption (CI)} an extension $\tilde{\mathcal{B}}$ is very good if and only if every martingale on $\mathcal{B}$ is a martingale on $\tilde{\mathcal{B}}$.
\end{lemma}
\noindent
The proof of this Lemma is similar to its one-parameter counterpart \cite[Lemma II.7.3]{JacodShiryaev}.

\subsection{Malliavin calculus for two-parameter Brownian motion}
\label{subsection:Malliavin calculus for two-parameter Brownian motion}

Let $(W_{(s,t)})_{(s,t) \in [0,1]^2}$ be a two-parameter Brownian motion defined on a probability basis $\left(\Omega,\mathcal{F},\left(\mathcal{F}_z\right)_{z\in[0,1]^2},\P\right)$. This process is a centered Gaussian process whose covariance given by,
$$\E\left[W_{(s,t)} W_{(s',t')}\right]=(s \wedge s') \, (t\wedge t'), \quad (s,t) \in [0,1]^2, \; (s',t')\in [0,1]^2.$$
The Malliavin calculus for general Gaussian processes has been described in \cite{Nualart3} and the reader can refer to it for a complete explanation about this topic. Here we give the definition of the Malliavin derivative and we present the integration by parts formula which is hardly used in Section \ref{section:Non-central limit theorem}.
\begin{definition}
Let $\mathcal{S}$ be the space of random variable $F$ of the form
\begin{equation}
\label{eq:cylindricalfunctional}
F=f(W(h_1),\ldots,W(h_n)),
\end{equation}
where $h_i$ is an element of $L^2([0,1]^2,dz)$ and $W(h_i)$ denotes the stochastic integral $W(h_i):=\int_{[0,1]^2} h_i(z) \; dz$ for $i=1,\cdots,n$ and $f:\real^n \to \real$ is infinitely continuously differentiable.\\
For $F$ of the form (\ref{eq:cylindricalfunctional}) we define the Malliavin derivative $DF$ of $F$ as the following $L^2([0,1]^2,dz)$-valued random variable,
\begin{equation}
\label{eq:MalliavinDerivative}
DF:=\sum_{i=1}^n \partial_i f(W(h_i)) \, h_i.
\end{equation} 
\end{definition}
\noindent
Here we give the Malliavin integration by parts formula (\cite[Lemma 1.2.1]{Nualart3}).  
\begin{lemma}
Let $F$ in $\mathcal{S}$ and $(u_{(s,t)})_{(s,t)\in[0,1]^2}$ in the domain of the divergence operator $\delta$. We have,
\begin{equation}
\label{eq:IntegrationByPartsFormula}
\E\left[F \delta(u) \right]=\E\left[\langle DF, \delta(u)\rangle_{L^2\left([0,1]^2,dz\right)}\right].
\end{equation}
$\delta$ is called the divergence operator (or the Skorohod integral of $u$) and it extends the It\^o integral since when $u$ is adapted to the filtration generated by the two-parameter Brownian motion $W$, 
$$\delta(u)=\int_{[0,1]^2} u_{(s,t)} \; dW_{(s,t)}.$$ 
\end{lemma}
\noindent Note also that the Malliavin derivative $D$ is a closable operator (see \cite[Proposition 1.2.1]{Nualart3}) and we denote by $Dom(D)$ its domain. Furthermore $D$ satisfy a chain rule property that is, for every $F$ and $G$ elements of $Dom(D)$ such that $F\,G$ belongs to $Dom(D)$ we have 
\begin{equation}
\label{eq:chainrule}
D(F G)=F\,DG + DF\,G.
\end{equation}
We end this section with the \textquotedblleft Heisenberg commutativity relationship\textquotedblright which enables the computation of the gradient of a It\^o stochastic integral, more precisely we have for a process $u$ such that the right hand side of (\ref{eq:Heisenbergcommutativityrelationship}) is well-defined (more details about assumption on $u$ can be found in \cite[(1.46)]{Nualart3}),
\begin{equation}
\label{eq:Heisenbergcommutativityrelationship}
D_{(s,t)}\delta(u) = u_{(s,t)} + \delta(D_{(s,t)} u), \quad (s,t)\in[0,1]^2. 
\end{equation}

\section*{Acknowledgement}
We thank the financial support of the \textit{European Social Fund (ESF)}.

\def\polhk#1{\setbox0=\hbox{#1}{\ooalign{\hidewidth
\lower1.5ex\hbox{`}\hidewidth\crcr\unhbox0}}} \def\cprime{$'$}

\end{document}